\documentstyle[12pt]{article}

\newcommand{\sect}[1]{\setcounter{equation}{0}\section{#1}}
\newcommand{\subsect}[1]{\subsection{#1}}

\renewcommand{\theequation}{\arabic{section}.\arabic{equation}}

\def\be{\begin{equation}}
\def\ee{\end{equation}}
\def\bea{\begin{eqnarray}}
\def\eea{\end{eqnarray}}

\def\osc{h_4}

\def\Osc{H_4}

\def\aa{N} 
\def\ap{{A_+}} 
\def\am{{A_-}}
\def\bb{M}

\def\aaa{n} 
\def\aap{{a_+}}
\def\aam{{a_-}}
\def\bbb{m}

\def\haaa{{\hat n}} 
\def\haap{{\hat{a}_+}}
\def\haam{{\hat{a}_-}}
\def\hbbb{\hat{m}}

\def\paaa{{n'}} 
\def\paap{{a'_+}}
\def\paam{{a'_-}}
\def\pbbb{{m'}}

\def\ppaaa{{n''}} 
\def\ppaap{{a''_+}}
\def\ppaam{{a''_-}}
\def\ppbbb{{m''}}

\def\xp{{\alpha_+}} 
\def\xm{{\alpha_-}}

\def\yp{{\beta_-}} 
\def\ym{{\beta_+}}

\def\xx{{\vartheta}} 

\def\yy{{\xi}}

\def\Ipp{{I_+}} 
\def\Imm{{I_-}} 
\def\II{{II}} 
 
\def\vv{{V}}

\def\a{\alpha} 
\def\bt{\tau}

\def\mm{\mu} 
\def\nn{\nu}

\def\>#1{{\bf #1}}                 
\def\1{\'{\i}}                           
\def\R{\rm I\kern-.2em R} 
\def\back{\!\!\!\!\!\!}

\begin{document}

\thispagestyle{empty}

\hfill{\ }
\ 
\vspace{2cm}

\begin{center} {\LARGE{\bf{Lie bialgebra quantizations of }}} 
  
 {\LARGE{\bf{ the oscillator
algebra and }}} 
  
  {\LARGE{\bf{ their universal $R$--matrices }}} 
 \end{center}

\bigskip\bigskip\bigskip

\begin{center} Angel Ballesteros and  Francisco J. Herranz
\end{center}

\begin{center} {\it {  Departamento de F\1sica, Universidad
de Burgos} \\   Pza. Misael Ba\~nuelos, \\
E-09001, Burgos, Spain}
\end{center}

\bigskip\bigskip\bigskip

\begin{abstract} 
All coboundary Lie bialgebras and their corresponding Poisson--Lie
structures are constructed for the oscillator algebra generated by 
$\{\aa,\ap,\am,\bb\}$. Quantum oscillator algebras are
derived from these bialgebras by using the Lyakhovsky and Mudrov
formalism and, for some cases, quantizations at both algebra
and group levels are obtained, including their universal $R$--matrices.
\end{abstract}

\newpage


\sect {Introduction}

Deformed Heisenberg and oscillator algebras  have recently focused many
investigations coming from different directions. Among them, we would
like to quote the construction of deformed statistics \cite{Greenberg},
the use of
$q$--Heisenberg algebras to describe composite particles \cite{Avancini}, the
description of certain classes of exactly solvable potentials in terms
of a $q$--Heisenberg dynamical symmetry \cite{Spiridonov}, the link
between deformed oscillator algebras and superintegrable systems
\cite{Bonatsosa,Bonatsosb} and the relations between these deformed
algebras and $q$--orthogonal polynomials
\cite{Floreanini}.

Quantum universal enveloping algebras (QUEA) are much more selective
deformations than general modifications of the commutation rules of a
given algebra. In particular, the interest of finding Hopf algebra
deformations of the oscillator algebra is twofold: firstly, because of 
the relevant role played by Hopf algebras to build up second
quantization, as it has been recently discussed in
\cite{Celeghini}. On the other hand, a quasitriangular quantum
oscillator algebra has been related to Yang--Baxter systems and link
invariants in
\cite{Sierra}.

The aim of this paper is to provide a systematic study of the quantum
universal enveloping oscillator algebras underlying possible further
generalizations of these results. A brief summary of the oscillator
algebra and group is given in section 2. Since  every QUEA defines
uniquely a Lie bialgebra structure on the undeformed algebra, in section
3 we obtain and classify all coboundary Lie bialgebra structures for the
harmonic oscillator algebra, as well as their corresponding  
Poisson--Lie brackets.  In section 4 we make use of the Lyakhovsky and
Mudrov formalism \cite{Lyak} in order to build up the deformed coproducts
linked to  all these Lie harmonic oscillator bialgebras.
A complete quantization (including universal
$R$--matrices) of two particular classes of
non-standard (triangular) bialgebras is provided: the former
is the natural ``extension" of the non-standard deformation of the 1+1
Poincar\'e algebra discussed in \cite{Tmatrix} and the latter is a new
three parameter quantization.

To our knowledge, the literature on Hopf algebra deformations of the
oscillator algebra includes only the deformation
given in \cite{Sierra,Celeghinidos} and some new results that have been
recently given in \cite{Vero} by computing the dual of an arbitrary
quantum oscillator group obtained by following an $R$--matrix approach
in a particular matrix representation (see \cite{FRT,BCGST,HLR}). Among
these known deformations, the former can be easily included within our
clasification at the Lie bialgebra level, and can thus be obtained
without making use of contraction procedures. On the other hand, our
method gives explicit (and universal) expressions for the oscillator QUE
algebras linked to the quantizations of \cite{Vero} which are
coboundaries. The procedure here outlined precludes cumbersome duality
computations and leads to rather simple candidates for universal
$R$--matrices.


\sect {Classical oscillator algebra and group} 

The oscillator Lie algebra $\osc$ is generated by
  $\{\aa,\ap,\am,\bb\}$   with Lie brackets
\be
[\aa,\ap]=\ap,\quad [\aa,\am]=-\am,\quad [\am,\ap]=\bb,\quad 
[\bb,\cdot\,]=0 .
\label{aa}
\ee
Besides the central generator $\bb$ there exists another Casimir
invariant:  \be
C=2\aa\bb - \ap\am -\am\ap .
\label{ab}
\ee 
A $3\times 3$ real matrix representation $D$ of (\ref{aa}) is given by:
\bea
&&D(\aa)=\left(\begin{array}{ccc}
 0 &0 & 0 \\ 0 & 1 & 0 \\ 0 & 0 & 0 
\end{array}\right),\qquad  
D(\ap)=\left(\begin{array}{ccc}
 0 &0 & 0 \\ 0 & 0 & 1 \\ 0 & 0 & 0 
\end{array}\right),\cr 
&&D(\am)=\left(\begin{array}{ccc}
 0 &1 & 0 \\ 0 & 0 & 0 \\ 0 & 0 & 0 
\end{array}\right),\qquad 
D(\bb)=\left(\begin{array}{ccc}
 0 &0 & 1 \\ 0 & 0 & 0 \\ 0 & 0 & 0 
\end{array}\right).
\label{ac}
\eea
The expression for a generic element of the oscillator group   $\Osc$
coming from this representation is: 
\bea
&&T^D=\exp\{\bbb D(\bb)\}\exp\{\aam D(\am)\}
\exp\{\aap D(\ap)\}\exp\{\aaa D(\aa)\}\cr
&&\quad =\left(\begin{array}{ccc}
 1 &\aam e^\aaa & \bbb+\aam\aap \\ 0 & e^\aaa & \aap \\ 0 & 0 & 1 
\end{array}\right) .
\label{ad}
\eea
The group law for the coordinates $\bbb$, $\aam$, $\aap$ and $\aaa$ is
obtained by means of matrix multiplication ${T^D}''={T^D}'\cdot {T^D}$:
\bea
&& \ppaaa=\aaa+\paaa ,\qquad \ppbbb=\bbb+\pbbb  
-\aam\paap  e^{-\paaa },\cr
&& \ppaap=\paap +\aap e^{\paaa },\qquad 
 \ppaam=\paam +\aam e^{-\paaa } .
\label{ae}
\eea
Left and right invariant vector fields are also deduced from
(\ref{ad}) and read: 
\be 
  X^L_\aa=\partial_\aaa,\quad X^L_\ap=e^\aaa\partial_\aap,\quad
  X^L_\am=e^{-\aaa}\partial_\aam - \aap e^{-\aaa}\partial_\bbb,\quad 
X^L_\bb=\partial_\bbb ;
\label{af}
\ee 
\be 
  X^R_\aa=\partial_\aaa + \aap \partial_\aap+ \aam \partial_\aam,\quad
   X^R_\ap=  \partial_\aap - \aam \partial_\bbb,\quad
  X^R_\am= \partial_\aam ,\quad 
X^R_\bb=\partial_\bbb .
\label{ag}
\ee 

The Heisenberg algebra can be seen as the subalgebra $\langle
\ap,\am,\bb\rangle$ of $\osc$ and the Heisenberg group $\langle
\aap,\aam,\bbb\rangle$ is recovered by taking the coordinate $\aaa\equiv
0$ in $\Osc$. Moreover, $\osc$ can be seen as a centrally extended (1+1)
Poincar\'e algebra (by
$\bb$). This fact will be useful in the
quantization process.


\sect {Coboundary oscillator Lie bialgebras}

Let $g$ be a Lie algebra and let $r$ be an element of $g\wedge g$. The
cocomutator $\delta:g\rightarrow g\wedge g$ given by 
\be
\delta(X)=[1\otimes X + X \otimes 1,\,  r],\quad 
X\in g, \label{bc}
\ee
defines a coboundary Lie
bialgebra $(g,\delta (r))$ if and only if $r$ fulfills the modified
classical Yang--Baxter equation (YBE)
\be
[X\otimes 1\otimes 1 + 1\otimes X\otimes 1 +
1\otimes 1\otimes X,[[r,r]]\, ]=0, \quad X\in g,
\label{mod}
\ee
where $[[r,r]]$ is the Schouten bracket defined by
\be
[[r,r]]:=[r_{12},r_{13}] + [r_{12},r_{23}] + 
[r_{13},r_{23}],
\label{bb}
\ee
and, if $r=r^{i j} X_i\otimes X_j$, we have denoted
$r_{12}=r^{i j} X_i\otimes X_j\otimes 1$,
$r_{13}=r^{i j} X_i\otimes 1\otimes X_j$ and
$r_{23}=r^{i j} 1\otimes X_i\otimes  X_j$.

When the $r$--matrix is such that $[[r,r]]=0$ (classical YBE), we shall
say that $(g,\delta (r))$ is a {\em non-standard} (or triangular) Lie
bialgebra. On the contrary, a solution $r$ of (\ref{mod}) with non
vanishing Schouten bracket will give rise to a so called  {\em standard}
Lie bialgebra. 

We recall that, if $g=Lie(G)$, the (unique) Poisson--Lie structure on
$C^\infty(G)$ linked to a fixed bialgebra $(g,\delta (r))$ is given by
the Sklyanin bracket 
\be
\{\Psi,\Phi\}=r^{\alpha\beta}\left(X_\alpha^L\Psi X_\beta^L\Phi
-X_\alpha^R\Psi X_\beta^R\Phi\right)\label{bd} , \qquad 
\Psi,\Phi\in C^\infty(G),
\ee
where $X_\alpha^L$ and $X_\bt^R$ are the left and right invariant vector
fields of $G$, respectively. 

In particular, for $\osc$ we shall consider an arbitrary element $r$,
which can be written in terms of six (real) coefficients:
\be 
 r=\xp\, \aa\wedge \ap+ \xm\, \aa\wedge \am + \xx\, \aa\wedge \bb
+ \yy\, \ap\wedge \am+ \ym\, \ap\wedge \bb 
 + \yp  \am\wedge \bb.
\label{ba}
\ee
It is a matter of computation to prove that the corresponding Schouten
bracket for $r$ (\ref{ba}) is 
\bea
&& \back [[r,r]]=
\xp\,(\yy+\xx)\, \aa\wedge\bb\wedge\ap +
\xm\,(\yy-\xx)\, \aa\wedge\bb\wedge\am \nonumber\\
&& \back  \qquad\quad -2\xp\xm\, \aa\wedge\ap\wedge\am +
(\xp\yp\, + \xm\ym- \yy^2)\,\bb\wedge\ap\wedge\am.
\label{bab}
\eea
From this expression follows that the modified
classical YBE (\ref{mod})
is fulfilled provided that:
\bea
&& \xp\xm=0,\cr
&& \xp(\yy+\xx)=0,\label{da}\\
&& \xm( \yy-\xx)=0 .\nonumber
\eea
The solutions of this system are splitted into three classes: 
$\xp\ne 0$, $\xm\ne 0$ and $\xp=\xm=0$. For each of them we
shall distinguish between non-standard ($[[r,r]]=0$) and standard
Lie bialgebras as follows.

\noindent
{\bf Type I$_+$.} If $\xp\ne 0$ we have  
$\xm=0$ and $\yy=-\xx$. The Schouten bracket reduces to:
\be
[[r,r]]=(\xp\yp-\xx^2)\, (\bb\wedge\ap\wedge\am).
\label{db}
\ee
Therefore if $\xp\yp-\xx^2\ne 0$ we have standard solutions and when
$\yp = \xx^2/\xp$ we are considering non-standard ones.

\noindent
{\bf Type I$_-$.} If $\xm\ne 0$ equations (\ref{da}) imply  
$\xp=0$ and $\yy=\xx$. The Schouten bracket is now:
\be
[[r,r]]=( \xm\ym-\xx^2)\, (\bb\wedge\ap\wedge\am).
\label{dbb}
\ee
Standard solutions are obtained when 
$\xm\ym-\xx^2\ne 0$, while non-standard ones correspond to 
$\ym = \xx^2/\xm$.

\noindent
{\bf Type II.} Finally, we  consider the case with 
$\xp=0$; if $\xm\ne 0$ we are again in type I$_-$, so we must take  also
$\xm=0$  in order to have three disjoint sets of solutions. In this case
equations (\ref{da}) are automatically satisfied and the 
  Schouten bracket is:
\be
[[r,r]]=-\yy^2\, \bb\wedge\ap\wedge\am.
\label{dbbb}
\ee
Then the condition $\yy\ne 0$ gives rise to standard solutions and
$\yy=0$ to non-standard ones.

All the information concerning this classification of coboundary
oscillator Lie bialgebras is summarized in table I. Poisson--Lie
structures for the oscillator group are deduced via the Sklyanin bracket
(\ref{bd}) and presented in table II.

Note that this classification is based in the use of skew-symmetric
$r$--matrices. This implies no loss of generality: given an arbitrary
element of $g\otimes g$, the map $\delta$ generated by (\ref{bc}) has to
be skew-symmetric to give rise to a Lie bialgebra. This amounts to impose
$Ad^{\otimes 2}$--invariance on the symmetric part of $r$ and, therefore,
$r$ will generate the same Lie bialgebra than its skew-symmetric part  
\cite{Tjin}. In particular, it can be easily checked that the more
general element $\eta$ of $\osc\otimes\osc$ such that 
\be
[X\otimes 1 +1\otimes X ,\eta]=0,\qquad X\in\{\aa,\ap,\am,\bb\} ;
\ee
is given by:
\be 
 \eta=\bt_1\, (\aa\otimes \bb +   \bb\otimes\aa - 
\ap\otimes \am - \am\otimes \ap) + \bt_2\,  \bb\otimes \bb,
\label{eta}
\ee
i.e., a linear combination of two terms directly related to the two
Casimirs of $\osc$.


\sect {Quantization}

In this section we first show how the Lyakhovsky and Mudrov (LM)
formalism \cite{Lyak} allows all the
cocommutators of the oscillator bialgebras previously found to generate
coassociative coproducts in a straightforward way. Afterwards, we shall 
construct commutation rules and universal quantum $R$--matrices for some
of these bialgebra quantizations.


\subsect {The Lyakhovsky--Mudrov formalism}

Let us start with a short resume of the LM formalism  which applies to an
associative algebra $E$ over $\bf  C$ with unit and generated by $n$
commuting elements $H_i$ and $m$ additional elements $X_j$. For any
$m\times m$ numerical matrix $\mu$, by $\mu\,H$ we understand the matrix
$\mu$ with all its entries multiplied by $H$. If $P$ is an $m\times m$
matrix with entries $p_{kl}\in E$, the $k$--th component of $P
\dot\otimes \vec X$ is defined as
\be 
(P\dot\otimes \vec X)_k=\sum_{l=1}^m p_{kl}\otimes X_l. 
\label{gf}
\ee
The main LM statement  \cite{Lyak}  is that $E$ can be endowed with a
coalgebra structure as follows (where we have denoted by $\sigma$ the
permutation map $\sigma (a\otimes b)=b\otimes a$):

\noindent
{\bf Proposition 1.} {\em Let $\{1, H_1,\dots,H_n,X_1,\dots,X_m\}$ a
basis of an associative algebra $E$ over $\bf  C$ verifying the conditions
\be
 [{H_i},{H_j}]=0,\qquad   i,j=1,\dots,n.
\label{ga}
\ee 
 Let $\mm_i$, $\nn_j$ $(i,j=1,\dots,n)$ be a set of
$m\times m$ complex matrices   such that
\be
[{\mm_i},{\nn_j}]=[{\mm_i},{\mm_j}]=[{\nn_i},{\nn_j}]=0 ,\qquad
  i,j=1,\dots,n.
\label{gb}
\ee
Let $\vec X$ be a   column vector with components $X_l$ $(l=1,\dots,m)$.
The coproduct and the counit
\bea
&&\Delta (1) =1 \otimes 1,\qquad 
\Delta (H_i) =1 \otimes H_i + H_i\otimes 1,\label{gc}\cr
&&\Delta (\vec X) = \exp (\sum_{i=1}^n{\mm_i H_i}) \dot \otimes \vec X +
\sigma\biggl( \exp (\sum_{i=1}^n{\nn_i H_i}) \dot\otimes \vec X \biggr),
\label{gd}\\
&&\epsilon(1)=1,\quad \epsilon(H_i)=\epsilon(X_l)=0,\qquad
i=1,\dots,n;\quad l=1,\dots m;
\label{ge}
\eea
endow $(E,\Delta,\epsilon)$ with a coalgebra structure.}

The resulting coalgebra can be seen as a multiparametric deformation
where the deformation parameters are the entries of the matrices $\mm_i$
and $\nn_j$. If we are able to find a compatible multiplication with
the coproduct (\ref{gd}) we will have finally obtained a quantum
algebra.

It is worth remarking that this formalism encodes in the set of matrices
$\mu_i$ and $\nu_j$ the whole coalgebra structure. In fact, the role
of these matrices is, essentially, to reflect the Lie bialgebra
underlying a given quantum deformation. This can be clearly appreciated
by taking the first order (in all the parameters) of (\ref{gd}):
\be
\Delta_{(1)} (\vec X) =   
(\sum_{i=1}^n{\mm_i H_i}) \dot \otimes \vec X +
\sigma\biggl(  (\sum_{i=1}^n{\nn_i H_i}) \dot\otimes \vec X \biggr)
\label{gg}
\ee
and recalling that the cocommutator $\delta$ corresponds to the
co-antisymmetric part of (\ref{gg}). It can be written in ``matrix" form
as:
\be \delta(\vec X)=\Delta_{(1)} (\vec X)-\sigma\circ\Delta_{(1)} (\vec
X). \label{gh}
\ee

We would like to emphasize the following points:

\noindent $\bullet$ The commuting elements $H_i$ are the primitive
generators. 

\noindent $\bullet$ The cocommutator $\delta(X_i)$ does not
contain terms of the form $H_i\wedge H_j$.

\noindent $\bullet$ The same cocommutator (\ref{gh}) can be obtained from
different choices of the matrices $\mu_i$ and $\nu_j$. This means that
  different sets of matrices might lead to right quantizations,   all of
them having the same first order terms in the deformation parameters.
Moreover, we can choose $\mu_i=0$ as a representative
of all these quantizations and we shall obtain  
\be
\delta  (\vec X) =  -  
  (\sum_{i=1}^n{\nn_i H_i}) \dot\wedge \vec X  
=-   (\sum_{i=1}^n{\nn_i H_i}) \dot\otimes \vec X +
\sigma\biggl(\sum_{i=1}^n{\nn_i H_i}) \dot\otimes \vec X\biggr).
\label{gi}
\ee

Now let us reverse somehow the LM formalism  trying to find in which
way the oscillator Lie bialgebras given in table I can be recovered by a
suitable choice of the matrices $\mu_i$ and $\nu_j$. Of course, the
benefit of such a situation is to be able to ``exponentiate" directly the
bialgebra (\ref{gi}) to a full coalgebra (\ref{gc}--\ref{ge}). 

Let us start with non-standard type I$_+$ oscillator bialgebras. By
denoting  $H_1\equiv \ap$, $H_2\equiv \bb$, $X_1\equiv \aa$, $X_2\equiv
\am$, we see that $[H_1,H_2]=0$ and there exists a term of the type
$H_1\wedge H_2\equiv\ap\wedge\bb$ within the cocommutator $\delta(\aa)$;
however, this obstruction can be circumvented by defining a new generator
in the form: 
\be
\aa'=\aa- (\ym/\xp)\bb .
\label{ggi}
\ee
 Hence, the cocommutators for the
non-primitive generators $\aa'$ and $\am$ can be written as
\be
 \delta\left(\begin{array}{c}
\aa' \\ \am  
\end{array}\right)=
\left(\begin{array}{cc}
-\xp\ap &0   \\ 0 & -\xp\ap 
\end{array}\right)\dot\wedge \left(\begin{array}{c}
\aa' \\ \am  
\end{array}\right)+
\left(\begin{array}{cc}
0 &(\xx^2/\xp)\bb   \\ -\xp\bb & -2\xx\bb 
\end{array}\right)\dot\wedge \left(\begin{array}{c}
\aa' \\ \am  
\end{array}\right)  
\label{gj}
\ee
In view of this expression, the matrices $\mm_i$ and $\nn_j$ can be
chosen as:
\be
\mm_1=\mm_2=\left(\begin{array}{cc}
0  &0   \\ 0 & 0 
\end{array}\right),\qquad
\nn_1=\left(\begin{array}{cc}
 \xp  &0   \\ 0 &  \xp 
\end{array}\right),\qquad \nn_2=\left(\begin{array}{cc}
0 &-\xx^2/\xp   \\  \xp  &  2\xx 
\end{array}\right).
\label{gk}
\ee
Now, the set of conditions of Proposition 1 are fulfilled, and we can use
this result to get  the coproducts:
\bea 
 &&\Delta\left(\begin{array}{c}
\aa' \\ \am  
\end{array}\right)=
\exp\left\{\left(\begin{array}{cc}
0 &0   \\ 0 & 0 
\end{array}\right)\right\}
\dot\otimes \left(\begin{array}{c}
\aa'\\ \am  
\end{array}\right) \cr
&&\qquad \qquad +
\sigma\left( \exp\left\{\left(\begin{array}{cc}
\xp\ap &- (\xx^2/\xp)\bb   \\ \xp\bb & \xp\ap+2\xx\bb 
\end{array}\right)\right\}\dot\otimes \left(\begin{array}{c}
\aa' \\ \am  
\end{array}\right)\right)\cr
&&
\!\!\!\!\!\!\!\!=\left(\begin{array}{c}
1\otimes \aa' + \aa'\otimes (1-\xx\bb)   e^{\xp\ap+\xx\bb}
-(\xx^2/\xp)\am\otimes \bb e^{\xp\ap+\xx\bb} \\ 
1\otimes \am + \am\otimes (1+\xx\bb)   e^{\xp\ap+\xx\bb}
+\xp \aa'\otimes \bb e^{\xp\ap+\xx\bb}
\end{array}\right)
\label{gm}
\eea
We can finally return to the initial basis elements, thus obtaining
a three-parameter QUEA (denoted by $U_{\xp,\xx,\ym}^{(\Ipp\, n)}(\osc)$)
such that:
\bea
&&
\Delta(\aa)=1\otimes \aa + \aa\otimes (1-\xx\bb)   e^{\xp\ap+\xx\bb}
-(\xx^2/\xp)\am\otimes \bb e^{\xp\ap+\xx\bb} \cr 
&&\qquad\qquad\qquad +(\ym/\xp)\bb\otimes (1-(1-\xx\bb)e^{\xp\ap+\xx\bb}),
\label{ggm}\\
&&
\Delta(\am)= 1\otimes \am + \am\otimes (1+\xx\bb)   e^{\xp\ap+\xx\bb}\cr
&&\qquad\qquad\qquad +(\xp \aa - \ym\bb)\otimes \bb e^{\xp\ap+\xx\bb}.
\nonumber
\eea

This quantization procedure can be applied to the remaining 
types of bialgebras in the same way. For the standard type I$_+$
bialgebras  we also use (\ref{ggi}), while for  the bialgebras of
type I$_-$ we introduce the new generator
\be
\aa'=\aa- (\yp/\xm)\bb .
\label{ggii}
\ee
On the contrary, no such a  kind of transformation is necessary to get
the coproducts for the Lie bialgebras of type II.

 The coproducts for the corresponding QUEA of the coboundary oscillator
Lie bialgebras of table I are written down in table III; we denote each
multiparametric quantum coalgebra by $U_{\a_i}^{(t\, m)}(\osc)$ where $t$ is the
type,  $m=s$ or $m=n$ according either to the standard or   non-standard
oscillator deformations and with  $\a_i$ being the  deformation
parameters.  The explicit expressions for the coproducts of
$U_{\xp,\xx,\ym,\yp}^{(\Ipp\, s)}(\osc)$ and 
$U_{\xm,\xx,\ym,\yp}^{(\Imm\, s)}(\osc)$ are rather complicated so we
keep their matrix forms written in terms of the generator $\aa'$ defined
by either (\ref{ggi}) or by (\ref{ggii}), respectively.

The final step in the quantization process of a fixed bialgebra is to
find the commutation relations compatible with its deformed coproduct
(counit and antipode can be obtained in the form explained in
\cite{Lyak}). In the following, we solve completely this problem and
construct the deformed Hopf algebras $U_{\a_i}^{(t\,m)}(\osc)$  for some
representative cases among the ones included in Table III.


\subsect {Non-standard type I$_+$: $U_{z}^{(n)}(\osc)$}

It is remarkable that the oscillator algebra with basis
$\{\aa,\ap,\am,\bb\}$ can be interpreted as an extended (1+1) Poincar\'e
algebra where $\aa$ is the boost  generator, $\ap$ and $\am$ generate
the translations along the light-cone and $\bb$ is the central
generator. This fact rises the question about whether it is possible to
implement in this extended case the universal (non-standard) quantum
deformation of the Poincar\'e algebra studied in \cite{Tmatrix} from a
$T$--matrix approach.

Let us consider  the non-standard oscillator bialgebras of  type I$_+$
with $\xx=\ym=0$ and $\xp\equiv z$. According to table I the Lie
bialgebra is characterized by commutation relations   (\ref{aa}),
classical $r$--matrix: \be 
 r=z\, \aa\wedge \ap ;
\label{ea}
\ee 
and cocommutators:
\bea
&&\delta(\ap)=0,\quad \delta(\bb)=0,\quad 
  \delta(\aa)=z\, \aa\wedge\ap  ,\cr
&& \delta(\am)=z\, (\am\wedge \ap+\aa\wedge\bb )   .
\label{eb}
\eea
Poisson--Lie brackets are easily deduced from table II:
\bea
&&\{\aaa,\aap\}=z\, (e^\aaa -1),\qquad
\{\aaa,\aam\}=0,\qquad
\{\aam,\aap\}=z\, \aam,\cr
&&\{\aaa,\bbb\}=z\, \aam,\qquad
\{\aap,\bbb\}=z\,\aam\aap ,\qquad
 \{\aam,\bbb\}=-z\, a_-^2  .
\label{ec}
\eea
A quantum deformation for this Lie bialgebra is given by the following
statement:

\noindent
{\bf Proposition 2.} {\em The coproduct $\Delta$, counit $\epsilon$,
antipode $\gamma$   
\bea
&&\Delta(\ap)=1\otimes \ap +\ap \otimes 1,\qquad 
 \Delta(\bb)=1\otimes \bb +\bb \otimes 1,\cr
&&\Delta(\aa)=1\otimes \aa +\aa \otimes e^{z\ap},\quad 
 \Delta(\am)=1\otimes \am +\am \otimes e^{z\ap}+z\aa\otimes \bb e^{z\ap} ;\cr
&& \label{ed}
\eea
\be
\epsilon(X)=0,\qquad X\in\{\aa,\ap,\am,\bb\};
\label{ee}
\ee
\bea
&&\gamma(\ap)=-\ap,\qquad \gamma(\bb)=-\bb,\cr 
&&\gamma(\aa)=-\aa  e^{-z\ap},\quad 
\gamma(\am)=-\am  e^{-z\ap}+ z \aa\bb  e^{-z\ap}; 
\label{ef}
\eea
and the commutation relations
\be 
 [\aa,\ap]=\frac{e^{z\ap}-1}{z},\quad [\aa,\am]=-\am
,\quad [\am,\ap]=\bb e^{z\ap},\quad 
[\bb,\cdot\,]=0 ,
\label{eg}
\ee 
determine a Hopf algebra (denoted by $U_z^{(n)}(\osc)$) which quantizes
the non-standard   bialgebra generated by the classical $r$--matrix
  (\ref{ea}).}

The coproduct  (\ref{ed}) is obtained from table III. Note that $\bb$
remains as a central generator. There is another element belonging to
the center of  $U_z^{(n)}\osc$ whose classical limit is (\ref{ab}),
namely \be
C_z=2\aa\bb+  
\frac{e^{-z\ap}-1}{z}\, \am +\am \, \frac{e^{-z\ap}-1}{z}.
\label{eh}
\ee

An important feature of the quantum algebra $U_z^{(n)}(\osc)$ is that
the generators   $\aa$ and $\ap$ form a Hopf subalgebra which coincides
exactly with the corresponding to the quantum Poincar\'e algebra of
\cite{Tmatrix}. We recall that for this Hopf subalgebra there is a
universal $R$--matrix given by: \be
R=\exp\{-z\ap\otimes\aa\}\exp\{z\aa\otimes\ap\} .
\label{ei}
\ee
Obviously,   (\ref{ei}) satisfies the quantum YBE for 
  $U_z^{(n)}\osc$, but moreover it verifies 
\be
\sigma\circ \Delta(X)=R\Delta(X)R^{-1},\qquad
\mbox{for}\ \ X\in\{\aa,\ap,\am,\bb\} .
\label{ej}
\ee
This assertion must be proved only for $\bb$ and $\am$; the proof for 
the former is trivial since it is a central generator, and for the latter
we have \bea
&&\exp\{z\aa\otimes\ap\}\Delta(\am)\exp\{-z\aa\otimes\ap\}
=1\otimes \am + \am \otimes 1=\Delta_0(\am) ,  \cr
&& \exp\{-z\ap\otimes\aa\}
\Delta_0(\am) 
\exp\{z\ap\otimes\aa\}
=\sigma \circ \Delta(\am). 
\eea

The fulfillment of relation  (\ref{ej}) allows to use the FRT approach
in order to get a quantum deformation  of $Fun(\Osc)$ by taking into
account that in the matrix representation (\ref{ac}) the universal
$R$--matrix   (\ref{ei}) collapses into: \be
D(R)=I\otimes I + z(D(\aa)\otimes D(\ap) - D(\ap)\otimes D(\aa)), 
\label{rr}
\ee
where $I$ is the $3\times 3$ identity matrix. Therefore, the Hopf
structure of the associated oscillator quantum group is given
by:  

\noindent
{\bf Proposition 3.} {\em The coproduct, counit,
antipode   
\bea
&&\Delta( \haaa)=1\otimes  \haaa +  \haaa \otimes 1,\cr
&&\Delta( \haap)=e^{ \haaa}\otimes  \haap + \haap \otimes 1,\cr
&&\Delta( \haam)=e^{- \haaa}\otimes  \haam + \haam \otimes 1,\cr
&&\Delta( \hbbb)=1\otimes  \hbbb + \hbbb \otimes 1
-e^{- \haaa} \haap\otimes  \haam ;
\label{ek}
\eea
\be
\epsilon(X)=0,\qquad X\in\{\haaa,\haap,\haam,\hbbb\};
\label{el}
\ee
\bea 
 &&\gamma(\haaa)=-\haaa,\qquad\quad \gamma(\haap)=-e^{-
\haaa}\haap,\nonumber\\
 &&\gamma(\haam)=-e^{\haaa}\haam ,\quad 
\gamma(\hbbb)=-\hbbb - (e^{-\haaa}\haap e^{\haaa})\haam ; 
\label{em}
\eea 
together with the commutation relations
\bea
&&[ \haaa, \haap]=z\, (e^{  \haaa} -1),\qquad
[ \haaa, \haam]=0,\qquad
[ \haam, \haap]=z\,  \haam,\cr
&&[ \haaa, \hbbb]=z\,  \haam,\qquad
[ \haap, \hbbb]=z\, \haam \haap ,\qquad
 [ \haam, \hbbb]=-z\,{ \haam}^2  ,
\label{en}
\eea
constitute a Hopf algebra denoted by  $Fun_z^{(n)}(\Osc)$.}

The coproduct (\ref{ek}), counit (\ref{el}) and antipode (\ref{em}) are
obtained from the relations $\Delta(T)=T\dot\otimes T$, $\epsilon(T)=I$
and $\gamma(T)=T^{-1}$, where $T\equiv T^D$ is  the generic element of the
oscillator group $\Osc$ (\ref{ad}). The commutation rules are deduced
from  $RT_1T_2=T_2T_1R$, where $T_1=T\otimes I$, $T_2=I\otimes T$ and $R$
given by (\ref{rr}).

The commutation relations (\ref{en})  can be seen as a Weyl
quantization $\{\, ,\, \}\to z^{-1}[\, ,\,]$ of the fundamental Poisson
brackets (\ref{ec}). It is also clear that the coalgebra structure of 
$Fun_z^{(n)}(\Osc)$ determined by the coproduct (\ref{ek}) and counit
(\ref{el}) is valid for any quantum group which deforms $Fun(\Osc)$.

Some features of this new quantum oscillator algebra can be emphasized:

\noindent $\bullet$ When 
the central extension $\bb$ and its corresponding quantum coordinate
$\hbbb$ vanish  all results concerning the quantum Poincar\'e algebra and
group given in \cite{Tmatrix} are recovered. In this sense,   the quantum
coordinates $\haaa$, $\haap$ and $\haam$ close a quantum Hopf subalgebra
which coincides exactly with the quantum Poincar\'e group just mentioned.

\noindent $\bullet$ The primitive generator involved in the
deformation is now $\ap$. This fact will be relevant at a representation
theory level and, consequently, from the point of view of the physical
properties of this deformed oscillator.

\noindent $\bullet$ The deformed Heisenberg
subalgebra generated by $\ap,\am$ and $\bb$ is not a Hopf subalgebra due
to the appearence of $\aa$ in $\Delta(\am)$. However, the Hopf
subalgebra structure can be recovered by working on a representation
where the central generator $\bb$ is expressed as a multiple of the
identity. In this situation, $\aa$ can be defined in terms of $\ap$ and
$\am$ by using the Casimir (\ref{eh}). In general, this type of
non-standard deformed bosons can be expected to build up $q$--boson
realizations of the already known non-standard quantum algebras
\cite{Ohn,Bey}.


\subsect {Non-standard type
II: $U_{\xx,\ym,\yp}^{(\II\,n)}(\osc)$}

The classical  $r$--matrix
\be
r=\xx\, \aa\wedge \bb + \ym\, \ap\wedge \bb + \yp\, \am\wedge \bb,
\label{ma}
\ee
originates a non-standard three-parametric oscillator bialgebra of type
II  whose cocommutators and associated Poisson--Lie brackets
appear  in tables I and  II, respectively. A quantum deformation of this
coboundary Lie bialgebra is given by: 

\noindent
{\bf Proposition 4.} {\em The Hopf algebra denoted by 
$U_{\xx,\ym,\yp}^{(\II\,n)}(\osc)$ which quantizes the oscillator
bialgebra generated by (\ref{ma}) has coproduct given in table III,
counit  (\ref{ee}), antipode
\bea
&&\gamma(\bb)=-\bb,\qquad \gamma(\ap)=-\ap  e^{\xx\bb},\qquad 
\gamma(\am)=-\am  e^{-\xx\bb},\cr
&&\gamma(\aa)=-\aa  - (\ym/\xx)\, \ap (1-e^{\xx\bb}) 
 - (\yp/\xx)\, \am (1-e^{-\xx\bb}) ; 
\label{mb}
\eea
and commutation relations 
\bea
&& [\aa,\ap]=\ap-\yp\, \vv(-\xx),\qquad  
[\aa,\am]=-\am-\ym \, \vv(\xx),\cr
&&[\am,\ap]=\bb,\qquad  [\bb,\cdot\,]=0 ,
\label{mc}
\eea 
where
\be
\vv(x):=\frac 1{x^2}(e^{x\bb}-1-x\bb) .
\label{md}
\ee}

Note that $\lim_{x\to 0}\vv(x)= \bb^2/2$. The quantum analogue of
(\ref{ab}): 
\be
C_{\xx,\ym,\yp}=2\aa\bb - \ap\am - \am\ap
+2\yp\,  \vv(-\xx) \am - 2 \ym\, \vv(\xx)\ap ,
\label{me}
\ee
belongs to the center of $U_{\xx,\ym,\yp}^{(\II\,n)}(\osc)$.

It is worth remarking that this quantum oscillator algebra can be
related with the results of \cite{Vero}:  
$U_{\xx,\ym,\yp}^{(\II\,n)}(\osc)$ can be seen as a Type II case with
$p\equiv \xx$, $q\equiv -\xx$, $b\equiv \yp$ and $c\equiv -\ym$.
Moreover,

\noindent
{\bf Proposition 5.} {\em The element
\bea
&&  R=\exp\{r\}
 =\exp\{\xx\, \aa\wedge \bb + \ym\, \ap\wedge \bb + \yp\,
\am\wedge \bb\} \cr
&& =\exp\{-\bb\otimes (\xx\, \aa  + \ym\, \ap + \yp\,
\am)\} \exp\{ (\xx\, \aa  + \ym\, \ap + \yp\,
\am) \otimes  \bb\} \cr
&& 
\label{mf}
\eea
satisfies both the quantum YBE and relation (\ref{ej}), so it is a
universal
$R$--matrix for $U_{\xx,\ym,\yp}^{(\II\,n)}(\osc)$.}

Since $\bb$ is a central generator, it is clear that (\ref{mf}) is a
solution of the quantum YBE. The proof for  property (\ref{ej}) is 
sketched  in Appendix A. In the   matrix representation (\ref{ac}) we get: 
\be
D(R)=I\otimes I + \xx\, D(\aa)\wedge D(\bb) + \ym\, D(\ap)\wedge D(\bb)
+ \yp\, D(\am)\wedge D(\bb).
\label{mg}
\ee
The FRT prescription leads now to another multiparametric quantum
deformation of the algebra of the smooth functions on the oscillator
group $Fun_{\xx,\ym,\yp}^{(\II\,n)}(\Osc)$, given by coproduct (\ref{ek}), counit (\ref{el}), antipode
(\ref{em}) and the non vanishing commutation rules 
\be 
[\haap,\hbbb]=-\xx\, \haap +\ym\, (e^\haaa -1) ,\qquad
[\haam,\hbbb]=\xx\, \haam +\yp\, (e^{-\haaa} -1).
\label{mh}
\ee 
The classical limit
(in the three parameters) is $Fun(\Osc)$ and, once more,
commutators (\ref{mh}) are   a Weyl quantization of the Poisson--Lie
brackets  written in table II.


\subsect {Standard type
II: $U_{z}^{(s)}(\osc)$}

The classical $r$--matrix which solves the classical YBE and underlies
the quantum oscillator  algebra obtained in \cite{Sierra,Celeghinidos} by
a contraction method can be expressed  in our notation as:
\be
r= - z\,(\aa\otimes\bb + \bb\otimes\aa) + 2z \,\am\otimes \ap.
\label{pa}
\ee
Its symmetric ($r_+$) and skew-symmetric ($r_-$) parts are:
\bea
&&\back\back\back\back r_+=(r+\sigma\circ r)/2=
z \,(\am\otimes \ap+\ap\otimes \am) - z\,(\aa\otimes\bb + \bb\otimes\aa),
\label{pb}\\
&&\back\back\back\back r_-=(r-\sigma\circ r)/2=
z \, \am\wedge \ap  .
\label{pc}
\eea
The symmetric part $r_+$ corresponds to the element $\eta$ (\ref{eta})
with the parameters $\bt_1=-z$ and $\bt_2=0$. On
the other hand, $r_-$ can be identified with a
standard classical $r$--matrix of type II with parameters
$\xx=\ym=\yp=0$ and $\yy\equiv -z$ (see table I). Both
the standard $r$--matrix (which coincides with $r_-$ (\ref{pc})) and the
non antisymmetric one (\ref{pa}) give rise to the same oscillator
bialgebra with cocommutators:  \be
\delta(\aa)=\delta(\bb)=0,\qquad \delta(\ap)=z\, \ap\wedge\bb,\qquad
\delta(\am)=z\, \am\wedge\bb .
\label{pd}
\ee
The associated non vanishing Poisson--Lie brackets (see table II) are:
\be
\{\aap,\bbb\}=z\, \aap,\qquad \{\aam,\bbb\}=z\, \aam .
\label{pe}
\ee

The quantum deformation of this coboundary oscillator bialgebra is
  given by:
 
\noindent
{\bf Proposition 6.} {\em The quantum algebra which quantizes the
standard bialgebra generated by (\ref{pa}) has a Hopf structure denoted
by 
$U_{z}^{(s)}(\osc)$ and characterized by the coproduct, counit, antipode  
\bea
&&\Delta(\aa)=1\otimes \aa +\aa \otimes 1,\quad 
\Delta(A_+')=e^{-z\bb}\otimes A_+' +A_+' \otimes 1,\cr
 &&
\Delta(\bb)=1\otimes \bb +\bb \otimes 1 
,\quad 
 \Delta(\am)=1\otimes \am +\am \otimes e^{z\bb}  ; 
  \label{pf}
\eea
\be
\epsilon(X)=0,\qquad X\in\{\aa,A_+',\am,\bb\};
\label{pg}
\ee
\be 
 \gamma(\aa)=-\aa,\quad \gamma(\bb)=-\bb,\quad  
 \gamma(A_+')=-A_+'  e^{z\bb},\quad 
\gamma(\am)=-\am  e^{-z\bb} ; 
\label{ph}
\ee 
together with the commutation relations
\be 
 [\aa,A_+']=A_+',\quad [\aa,\am]=-\am
,\quad [\am,A_+']=\frac{\sinh(z\bb)}{z},\quad 
[\bb,\cdot\,]=0 .
\label{pi}
\ee }

The quantum Casimir is:
\be
C_z=2\aa\, \frac{\sinh(z\bb)}{z} - A_+'\am -\am A_+' .
\label{pj}
\ee

The coproducts (\ref{pf}) are just those given in table III but written
in terms of a new generator
$A_+'=e^{\yy\bb}\ap$ where $\yy =- z$.
In this case the universal $R$--matrix adopts a much
 simpler form than the one  already known from
\cite{Sierra,Celeghinidos}. Namely,
\bea
&&  R =\exp\{-z\,(\aa\otimes\bb + \bb\otimes\aa)\}\exp\{2z \,\am\otimes
A_+'\}\cr
&&\quad\quad = \exp\{-z\, \aa\otimes\bb\}\exp\{-z\,
\bb\otimes\aa\}\exp\{2z \,\am\otimes A_+'\}. \label{pk}
\eea
It is worth remarking that all the 
quantum $R$--matrices given in this section are obtained via a 
straightforward exponentiation process from their classical
counterparts (compare, for instance, (\ref{pk}) to (\ref{pa})).

The FRT prescription can be applied leading to the commutation rules of
the quantum group $Fun_z^{(s)}(\Osc)$ by taking into account that
(\ref{pk}) in the matrix representation (\ref{ac}) is just 
\be
D(R)=I\otimes I +  2z \,D(\am)\otimes D(A_+') 
- z\,(D(\aa)\otimes D(\bb) + D(\bb)\otimes D(\aa)) ;
\label{pl}
\ee
(note that $D(A_+')\equiv D(\aa)$). In this way, the non vanishing
commutators of
$Fun_z^{(s)}(\Osc)$ read:
\be
[\haap,\hbbb]=z\, \haap,\qquad [\haam,\hbbb]=z\, \haam ,
\label{pm}
\ee
and correspond to a Weyl quantization of the Poisson--Lie brackets
(\ref{pe}).


\sect {Concluding remarks}

We have presented a systematic procedure in order to study the
coboundary Lie bialgebras of the oscillator algebra. The first order
deformations given by the corresponding cocomutators have been used to
construct, by a sort of ``exponentiation" process, multiparametric
quantum deformations of the oscillator algebra. We point out that we have
not treated the question of the equivalence of the coboundary oscillator
bialgebras we have obtained, indeed this is actually a problem by itself.
For instance, from an algebraic point of view, bialgebras of types I$_+$
and I$_-$ can be related by interchanging generators $\ap$ and $\am$,
although this result is not so straightforward if we look at their
corresponding Poisson--Lie groups.

It is worth stressing that, in the case here analysed, the complete
(and rich) classification of the classical $r$--matrices (and, therefore,
of the corresponding Poisson structures on the oscillator group) is
easily obtained. This seems to indicate that, at least for Lie algebras
with a   low enough dimension, the complete solution of the modified
classical YBE for an arbitrary skew element of $g\otimes g$ can be
explicitly deduced      
  giving rise to a great amount of new results.

 This kind of procedure is complementary (and dual) to
that developed in \cite{Vero}, since it allows us to focus on
the deformation at the quantum algebra level and looking for universal
quantum $R$--matrices. In fact, given a skew solution $r$ of the 
modified classical YBE and a matrix representation $D$ of the quantum
algebra, the element $D(R)=1+z\,D(r)$ will lead us to the corresponding
$R$--matrix method. 

This approach can be seen
as a part of a research program that, in order to construct and study
quantum algebras, tries to extract as much information as possible from
the associated Lie bialgebras (as far as contraction methods are
concerned, see for instance \cite{LBC}). It would be interesting to
apply    it to other physically interesting algebras whose coboundary
bialgebra structures are not well known, among them, we would like to
mention the Schr\"odinger, optical and Galilean algebras, also with the
aim of obtaining some (universal) quantum deformations.



\bigskip
\bigskip

\noindent
{\large{{\bf Acknowledgements}}}

\bigskip

The authors acknowledge   M. Santander and   M.A. del Olmo a careful
reading of the manuscript and many helpful suggestions as well as the
referees for some pertinent comments. This work has been partially
supported by DGICYT (Projects  PB92--0255 and PB94--1115) from the
Ministerio de Educaci\'on y Ciencia de Espa\~na.

\bigskip
\bigskip



\noindent
{\large{{\bf Appendix A}}}

\appendix

\setcounter{equation}{0}

\renewcommand{\theequation}{A.\arabic{equation}}

\bigskip

\noindent
The main steps necessary to prove that the $R$--matrix (\ref{mf})
verifies the property (\ref{ej}) for the generators $\ap$, $\am$ and
$\aa$ (for
$\bb$ the proof is trivial) are as follows. We perform the computations
by writing  the $R$--matrix in terms of two exponentials
$R=\exp\{-\bb \otimes W\}\exp\{W\otimes \bb\}$, where $W\equiv  \xx  \aa 
 + \ym  \ap + \yp \am$. We note that
\bea
&&\exp\{W\otimes  \bb\}\Delta(\ap)\exp\{- W\otimes  \bb\}\cr
&&=1\otimes \ap+ \ap\otimes 1 -\yp \vv(-\xx)\otimes (1-e^{-\xx\bb}) 
  +  (\yp/\xx) \bb\otimes (1-e^{-\xx\bb})\cr
&&= \Delta_0(\ap) + (\yp/\xx^2)\, (1-e^{-\xx\bb}) \otimes
(1-e^{-\xx\bb}).\label{xa}
\eea
Since the second term of (\ref{xa}) is central, we compute:
\bea
&&\exp\{- \bb\otimes  W \}\Delta_0(\ap)\exp\{\bb\otimes  W \}\cr
&&=e^{-\xx\bb}\otimes \ap+ \ap\otimes 1
 +\yp  (1-e^{-\xx\bb}) \otimes \vv(-\xx)
  -  (\yp/\xx)  (1-e^{-\xx\bb})\otimes \bb\cr
&&= \sigma \circ\Delta(\ap) - (\yp/\xx^2) \, (1-e^{-\xx\bb}) \otimes
(1-e^{-\xx\bb}).\label{xb}
\eea
From these expressions  $R\,\Delta(\ap)\,R^{-1}=
\sigma \circ\Delta(\ap)$ is easily derived. The proof for $\am$ is
rather similar, and for the generator
$\aa$ we shall have:
\bea
&&\exp\{W\otimes  \bb\}\Delta(\aa)\exp\{- W\otimes  \bb\}\cr
&&=1\otimes \aa+ \aa\otimes 1 
-(\ym\yp/\xx^2) \left\{ \xx \vv(\xx)\otimes (1-e^{\xx\bb}) 
+\bb\otimes (1-e^{-\xx\bb})\right\}\cr
&&\qquad 
+(\ym\yp/\xx^2) \left\{ \xx \vv(-\xx)\otimes (1-e^{\xx\bb}) 
-\bb\otimes (1-e^{-\xx\bb})\right\} \cr
&&=\Delta_0(\aa) 
+(\ym\yp/\xx^3) \left\{(1-e^{\xx\bb}) \otimes
(1-e^{\xx\bb}) -(1-e^{-\xx\bb}) \otimes
(1-e^{-\xx\bb})\right\} .\cr
&&
\label{xc}
\eea
Now we compute
\bea
&&\exp\{- \bb\otimes  W \}\Delta_0(\aa)\exp\{\bb\otimes  W \}\cr
&&=\Delta_0(\aa)+(\ym/\xx)\, (1-e^{-\xx\bb})\otimes \ap
+(\yp/\xx)\, (1-e^{\xx\bb})\otimes \am\cr
&&\quad +(\ym\yp/\xx ) \left\{(1-e^{\xx\bb}) 
\otimes    \vv(\xx) - (1-e^{-\xx\bb}) \otimes   \vv(-\xx)\right\}\cr
&&\quad +(\ym\yp/\xx^2)\, 
 (2-e^{\xx\bb}-e^{-\xx\bb})\otimes \bb \cr
&&=\sigma\circ\Delta(\aa)
-(\ym\yp/\xx^3) \left\{(1-e^{\xx\bb}) \otimes
(1-e^{\xx\bb}) -(1-e^{-\xx\bb}) \otimes
(1-e^{-\xx\bb})\right\},  \cr
&& \label{xd}
\eea
to obtain again  $R\,\Delta(\aa)\,R^{-1}=
\sigma \circ\Delta(\aa)$.

\bigskip
\bigskip

\newpage

\medskip

{\footnotesize

 \noindent
{{\bf Table I.} Coboundary oscillator Lie bialgebras.}
\smallskip

\noindent
\begin{tabular}{|c|l|l|}
\hline
 &\multicolumn{2}{c|}{\bf Type I$_+$\ \ \ \ \ }\\
\cline{2-3}
 & \multicolumn{1}{c|}{\em Standard}&
\multicolumn{1}{c|}{\em Non-standard}\\
&\multicolumn{1}{c|}{($\xp\ne0$ and $ \xp\yp-\xx^2\ne 0$)}
&\multicolumn{1}{c|}{($\xp\ne0$)}\\
\hline
$r$&\multicolumn{1}{l|}{$
 \xp\, \aa\wedge \ap+  \xx\, (\aa\wedge \bb - \ap\wedge \am)$}
&\multicolumn{1}{l|}{$  \xp\, \aa\wedge \ap+ 
 \xx\, (\aa\wedge \bb - \ap\wedge \am)$}\\
&\multicolumn{1}{l|}{$
\qquad  + \ym\, \ap\wedge \bb 
 + \yp \, \am\wedge \bb$}
&\multicolumn{1}{l|}{$\qquad 
 + \ym\, \ap\wedge \bb 
 + (\xx^2/\xp) \, \am\wedge \bb$}\\
\hline
$\delta(\aa)$& $\xp\, \aa\wedge\ap  -
\yp\, \am\wedge \bb +\ym\, \ap\wedge\bb$&
$\xp\, \aa\wedge\ap  -
(\xx^2/\xp)\, \am\wedge \bb +\ym\, \ap\wedge\bb$\\
$\delta(\ap)$&$\quad 0$&$\quad 0$\\
$\delta(\am)$&$\xp\, (  \am\wedge \ap +\aa\wedge\bb) 
+2\xx\, \am\wedge\bb$&
$\xp\, ( \am\wedge \ap +\aa\wedge\bb ) +2\xx\, \am\wedge\bb$\\
$\delta(\bb)$&$\quad 0$&$\quad 0$\\
\hline
\hline
 &\multicolumn{2}{c|}{\bf Type I$_-$\ \ \ \ \ }\\
\cline{2-3}
 & \multicolumn{1}{c|}{\em Standard}&
\multicolumn{1}{c|}{\em Non-standard}\\
&\multicolumn{1}{c|}{($\xm\ne0$ and $  \xm\ym-\xx^2\ne 0$)}
&\multicolumn{1}{c|}{($\xm\ne0$)}\\
\hline
$r$&\multicolumn{1}{l|}{$
  \xm\, \aa\wedge \am+  \xx\, (\aa\wedge \bb + \ap\wedge \am)$}
&\multicolumn{1}{l|}{$ \xm\, \aa\wedge \am+
\xx\, (\aa\wedge \bb + \ap\wedge \am) $}\\
&\multicolumn{1}{l|}{$
\qquad   + \ym\, \ap\wedge \bb 
 +\yp \, \am\wedge \bb$}
&\multicolumn{1}{l|}{$\qquad  + (\xx^2/\xm)\, \ap\wedge \bb 
 +\yp \, \am\wedge \bb$}\\
\hline
$\delta(\aa)$& $-\xm\, \aa\wedge\am + \ym\, \ap\wedge\bb -
\yp\, \am\wedge \bb$&
$-\xm\, \aa\wedge\am + (\xx^2/\xm)\, \ap\wedge\bb -
\yp\, \am\wedge \bb$\\
$\delta(\ap)$&$-\xm\, ( \ap\wedge \am + \aa\wedge\bb) 
-2\xx\, \ap\wedge\bb$&$-\xm\, (\ap\wedge \am +\aa\wedge\bb ) 
-2\xx\, \ap\wedge\bb$\\
$\delta(\am)$&$\quad 0$&
$\quad 0$\\
$\delta(\bb)$&$\quad 0$&$\quad 0$\\
\hline
\hline
 &\multicolumn{2}{c|}{\bf Type II\ \ \ \ \ }\\
\cline{2-3}
 & \multicolumn{1}{c|}{\em Standard}&
\multicolumn{1}{c|}{\em Non-standard}\\
&\multicolumn{1}{c|}{($\yy\ne0$ )}
&\multicolumn{1}{c|}{ }\\
\hline
$r$&\multicolumn{1}{l|}{$
     \xx\, \aa\wedge \bb  + \yy\, \ap\wedge \am$}
&\multicolumn{1}{l|}{$    \xx\, \aa\wedge \bb  
 + \ym\, \ap\wedge \bb 
 +\yp \, \am\wedge \bb$}\\
&\multicolumn{1}{l|}{$
\qquad    
 + \ym\, \ap\wedge \bb 
 +\yp \, \am\wedge \bb$}
&\multicolumn{1}{l|}{$ $}\\
\hline
$\delta(\aa)$& $\ym\, \ap\wedge\bb -\yp\, \am\wedge\bb$&
$\ym\, \ap\wedge\bb -\yp\, \am\wedge\bb$\\
$\delta(\ap)$&$-(\xx+\yy)\, \ap\wedge\bb$&
$-\xx\, \ap\wedge\bb$\\
$\delta(\am)$&$(\xx-\yy)\, \am\wedge\bb$&
$\xx\,
\am\wedge\bb$\\
$\delta(\bb)$&$\quad 0$&$\quad 0$\\
\hline
\end{tabular}}

 \newpage

{\footnotesize

 \noindent
{{\bf Table II.} Poisson--Lie brackets on the  oscillator group.}
\smallskip

\noindent
\begin{tabular}{|c|l|l|}
\hline
 &\multicolumn{2}{c|}{\bf   Type I$_+$\ \ \ \ \ }\\
\cline{2-3}
 & \multicolumn{1}{c|}{\em Standard}&
\multicolumn{1}{c|}{\em Non-standard}\\
&\multicolumn{1}{c|}{($\xp\ne0$ and $\xp\yp-\xx^2\ne 0$)}
&\multicolumn{1}{c|}{($\xp\ne0$)}\\
\hline
$\{\aaa,\aap\}$& $\xp\, (e^\aaa -1)$&$\xp\, (e^\aaa -1)$\\
$\{\aaa,\aam\}$&$\quad 0$&$\quad 0$\\
$\{\aam,\aap\}$&$\xp\, \aam$&
$\xp\, \aam$\\
$\{\aaa,\bbb\}$&$\xp\, \aam$&$\xp\, \aam$\\
$\{\aap,\bbb\}$&$\xp\,\aam\aap +\ym\, (e^\aaa -1)$&
$\xp\,\aam\aap +\ym\, (e^\aaa -1)$\\
$\{\aam,\bbb\}$&$-\xp\, a_-^2 +2\xx\, \aam +
\yp\,  (e^{-\aaa} -1)$&$-\xp\, a_-^2 +2\xx\, \aam +
(\xx^2/\xp)\,  (e^{-\aaa} -1)$\\
\hline
\hline
 &\multicolumn{2}{c|}{  \bf Type I$_-$\ \ \ \ \ }\\
\cline{2-3}
 & \multicolumn{1}{c|}{\em Standard}&
\multicolumn{1}{c|}{\em Non-standard}\\
&\multicolumn{1}{c|}{($\xm\ne0$ and $  \xm\ym-\xx^2\ne 0$)}
&\multicolumn{1}{c|}{($\xm\ne0$)}\\
\hline
$\{\aaa,\aap\}$& $\quad 0$&
$\quad 0$\\
$\{\aaa,\aam\}$&$\xm\, (e^{-\aaa} -1)$&$\xm\, (e^{-\aaa} -1)$\\
$\{\aam,\aap\}$&$\xm\, \aap$&
$\xm\, \aap$\\
$\{\aaa,\bbb\}$&$-\xm\, \aap e^{-\aaa}$&$-\xm\, \aap e^{-\aaa}$\\
$\{\aap,\bbb\}$&$-2\xx\, \aap +\ym\, (e^\aaa -1)$&
$-2\xx\, \aap +(\xx^2/\xm)\, (e^\aaa -1)$\\
$\{\aam,\bbb\}$&$\yp\,  (e^{-\aaa} -1)$&$\yp\,  (e^{-\aaa} -1)$\\
\hline
\hline
 &\multicolumn{2}{c|}{  \bf Type II\ \ \ \ \ }\\
\cline{2-3}
 & \multicolumn{1}{c|}{\em Standard}&
\multicolumn{1}{c|}{\em  Non-standard}\\
&\multicolumn{1}{c|}{($\yy\ne0$ )}
&\multicolumn{1}{c|}{ }\\
\hline
$\{\aaa,\aap\}$& $\quad 0$&
$\quad 0$\\
$\{\aaa,\aam\}$&$\quad 0$&$\quad 0$\\
$\{\aam,\aap\}$&$\quad 0$&
$\quad 0$\\
$\{\aaa,\bbb\}$&$\quad 0$&$\quad 0$\\
$\{\aap,\bbb\}$&$-(\xx+\yy)\, \aap +\ym\, (e^\aaa -1)$&
$-\xx\, \aap +\ym\, (e^\aaa -1)$\\
$\{\aam,\bbb\}$&$(\xx-\yy)\, \aam +\yp\, (e^{-\aaa} -1)$&$
\xx\, \aam +\yp\, (e^{-\aaa} -1)$\\
\hline
\end{tabular}}

\bigskip

\newpage

\bigskip

{\footnotesize

\noindent
{{\bf Table III.} Coproducts for QUEA of the oscillator algebra.}
\smallskip

\noindent
\begin{tabular}{|ll|}
\hline
 \multicolumn{2}{|c|}{\bf Type I$_+$}\\
\hline
\qquad\qquad {\em Standard}:\quad 
$U_{\xp,\xx,\ym,\yp}^{(\Ipp\, s)}(\osc)$\qquad 
&\qquad($\xp\ne0$ and $ \xp\yp-\xx^2\ne 0$)\\
\hline
 \multicolumn{2}{|l|} {$\Delta(\ap)=1\otimes \ap + \ap \otimes 1\qquad 
\Delta(\bb)=1\otimes \bb + \bb \otimes 1$}\\
 \multicolumn{2}{|l|}{$\Delta\left(\begin{array}{c}
\aa' \\ \am  
\end{array}\right)=
\left(\begin{array}{c}
1\otimes \aa'  \\ 1\otimes \am
\end{array}\right) +
\sigma\left( \exp\left\{\left(\begin{array}{cc}
\xp\ap &- \yp \bb   \\ \xp\bb & \xp\ap+2\xx\bb 
\end{array}\right)\right\}\dot\otimes \left(\begin{array}{c}
\aa' \\ \am  
\end{array}\right)\right)$}\\
\multicolumn{2}{|l|}{$\aa'=\aa-(\ym/\xp)\bb$}\\
 \hline
\qquad\qquad {\em Non-standard}:\quad 
$U_{\xp,\xx,\ym}^{(\Ipp\,n)}(\osc)$\qquad &\qquad($\xp\ne0$)\\
\hline
 \multicolumn{2}{|l|}{$\Delta(\ap)=1\otimes \ap + \ap \otimes 1\qquad 
\Delta(\bb)=1\otimes \bb + \bb \otimes 1$}\\
 \multicolumn{2}{|l|}{$\Delta(\aa)=1\otimes \aa + \aa\otimes (1-\xx\bb) 
\exp\{\xp\ap+\xx\bb\}
-(\xx^2/\xp)\am\otimes \bb \exp\{\xp\ap+\xx\bb\} $}\\ 
 \multicolumn{2}{|l|}{$\qquad\qquad +(\ym/\xp)\bb\otimes 
(1-(1-\xx\bb)  \exp\{\xp\ap+\xx\bb\} )$}\\
 \multicolumn{2}{|l|}{$\Delta(\am)= 1\otimes \am + \am\otimes (1+\xx\bb)  
\exp\{\xp\ap+\xx\bb\}$}\\
 \multicolumn{2}{|l|}{$\qquad\qquad 
+(\xp \aa-\ym\bb)\otimes \bb \exp\{\xp\ap+\xx\bb\}$}\\
\hline
\hline
 \multicolumn{2}{|c|}{\bf Type I$_-$}\\
\hline
\qquad\qquad {\em Standard}:\quad 
$U_{\xm,\xx,\ym,\yp}^{(\Imm\,s)}(\osc)$\qquad 
&\qquad($\xm\ne0$ and $ \xm\ym-\xx^2\ne 0$)\\
\hline
 \multicolumn{2}{|l|} {$\Delta(\am)=1\otimes \am + \am \otimes 1\qquad 
\Delta(\bb)=1\otimes \bb + \bb \otimes 1$}\\
 \multicolumn{2}{|l|}{$\Delta\left(\begin{array}{c}
\aa' \\ \ap  
\end{array}\right)=
\left(\begin{array}{c}
1\otimes \aa'  \\ 1\otimes \ap
\end{array}\right) +
\sigma\left( \exp\left\{\left(\begin{array}{cc}
-\xm\am &  \ym \bb   \\ -\xm\bb & -\xm\am-2\xx\bb 
\end{array}\right)\right\}\dot\otimes \left(\begin{array}{c}
\aa' \\ \ap  
\end{array}\right)\right)$}\\
 \multicolumn{2}{|l|} {$\aa'=\aa-(\yp/\xm)\bb $}\\
\hline
\qquad\qquad {\em Non-standard}:\quad 
$U_{\xm,\xx,\yp}^{(\Imm\,n)}(\osc)$\qquad &\qquad($\xm\ne0$)\\
\hline
 \multicolumn{2}{|l|}{$\Delta(\am)=1\otimes \am + \am \otimes 1\qquad 
\Delta(\bb)=1\otimes \bb + \bb \otimes 1$}\\
 \multicolumn{2}{|l|}{$\Delta(\aa)=1\otimes \aa + \aa\otimes (1+\xx\bb) 
\exp\{-\xm\am-\xx\bb\}
+(\xx^2/\xm)\ap\otimes \bb \exp\{-\xm\am-\xx\bb\} $}\\ 
 \multicolumn{2}{|l|}{$\qquad\qquad +(\yp/\xm)\bb\otimes 
(1-(1+\xx\bb)  \exp\{-\xm\am-\xx\bb\} )$}\\
 \multicolumn{2}{|l|}{$\Delta(\ap)= 1\otimes \ap + \ap\otimes (1-\xx\bb)  
\exp\{-\xm\am-\xx\bb\}$}\\
 \multicolumn{2}{|l|}{$\qquad\qquad
-(\xm \aa - \ym\bb)\otimes \bb \exp\{-\xm\am-\xx\bb\}$}\\
\hline
\hline
 \multicolumn{2}{|c|}{\bf Type II}\\
\hline
\qquad\qquad {\em Standard}:\quad 
$U_{\xx,\yy,\ym,\yp}^{(\II\, s)}(\osc)$\qquad 
&\qquad($\yy\ne0$)\\
\hline
 \multicolumn{2}{|l|} {$\Delta(\bb)=1\otimes \bb + \bb \otimes 1$}\\
 \multicolumn{2}{|l|} {$
\Delta(\ap)=1\otimes \ap + \ap \otimes \exp\{-(\xx+\yy)\bb\} $}\\
 \multicolumn{2}{|l|} 
{$\Delta(\am)=1\otimes \am + \am \otimes \exp\{(\xx-\yy)\bb\} $}\\
 \multicolumn{2}{|l|}{$\Delta(\aa)=1\otimes \aa + \aa\otimes 1
+ ({\ym}/({\xx+\yy})) \ap\otimes (1-\exp\{-(\xx+\yy)\bb\})$}\\ 
 \multicolumn{2}{|l|}{$
\qquad\qquad +({\yp}/({\xx-\yy})) \am\otimes
(1-\exp\{(\xx-\yy)\bb\})$}\\ 
\hline
\qquad\qquad {\em Non-standard}:\quad 
$U_{\xx,\ym,\yp}^{(II\,n)}(\osc)$\qquad 
&\qquad \\
\hline
 \multicolumn{2}{|l|} {$\Delta(\bb)=1\otimes \bb + \bb \otimes 1$}\\
 \multicolumn{2}{|l|} {$
\Delta(\ap)=1\otimes \ap + \ap \otimes \exp\{-\xx\bb\}$}\\
 \multicolumn{2}{|l|} {$
\Delta(\am)=1\otimes \am + \am \otimes \exp\{\xx\bb\} $}\\
 \multicolumn{2}{|l|}{$\Delta(\aa)=1\otimes \aa + \aa\otimes 1
+({\ym}/{\xx}) \ap\otimes (1-\exp\{-\xx\bb\})
+({\yp}/{\xx}) \am\otimes (1-\exp\{\xx\bb\})  $}\\ 
\hline
\end{tabular}}

\bigskip

\newpage


\end{document}